\documentclass[12pt]{article}
\usepackage{amsmath,amssymb,amsthm}
\usepackage{fullpage}

\newtheorem{thm}{\textbf{Theorem}}

\newtheorem{lemma}{Lemma}
\newcommand{\dis}{\displaystyle}

\begin{document}

\title{A Direct Proof of Landen's Transformation }
\author{Mark B. Villarino\\
Depto.\ de Matem\'atica, Universidad de Costa Rica,\\
2060 San Jos\'e, Costa Rica}
\date{\today}

\maketitle

 \begin{abstract}
 We prove Landen's transformation by direct elementary transformation of infinite series.
\end{abstract}
\section{Introduction}

Our paper contains a direct elementary proof of the following series transformation:

\begin{thm} (\textbf{Landen's Transformation})
If $0\leqslant x\leqslant 1$, then\begin{equation*}
\fbox{$\begin{array}{lll} 
    \dis(1+\sqrt{x})\left\{1-\frac{1}{2^{2}}\left(\frac{2\sqrt[4]{x}}{1+\sqrt{x}}\right)^{2}-\frac{1^{2}\cdot 3}{2^{2}\cdot 4^{2}}\left(\frac{2\sqrt[4]{x}}{1+\sqrt{x}}\right)^{4}-\frac{1^{2}\cdot 3^{2}\cdot 5}{2^{2}\cdot 4^{2}\cdot 6^{2}}\left(\frac{2\sqrt[4]{x}}{1+\sqrt{x}}\right)^{6}-\cdots\right\}=
    &\\
    &\\=\dis1+\left(\frac{1}{2}\right)^{2}x+\left(\frac{1\cdot 1}{2\cdot 4}\right)^{2}x^{2}+\left(\frac{1\cdot 1\cdot 3}{2\cdot 4\cdot 6}\right)^{2}x^{3}+\left(\frac{1\cdot 1\cdot 3\cdot 5}{2\cdot 4\cdot 6\cdot 8}\right)^{2}x^{4}+\cdots &  \\\end{array} $}
\end{equation*}\end{thm}\hfill$\Box$

The interest in this theorem is that it gives a direct proof the equality of the following two formulas for the perimeter of an ellipse with semi-axes $a$ and $b$.
\begin{thm} (\textbf{The Perimeter of an Ellipse}) The perimeter, $p$, of an ellipse with semi-major axes $a$ and $b$ and eccentricity $e$ is:\begin{equation*}
\fbox{$\begin{array}{lll} 
    \dis p=2a\pi\left\{1-\frac{1}{2^{2}}e^{2}-\frac{1^{2}\cdot 3}{2^{2}\cdot 4^{2}}e^{4}-\frac{1^{2}\cdot 3^{2}\cdot 5}{2^{2}\cdot 4^{2}\cdot 6^{2}}e^{6}-\cdots\right\}
    \\
    \\\ \ =\pi(a+b)\left\{\dis1+\left(\frac{1}{2}\right)^{2}\left(\frac{a-b}{a+b}\right)^{2}+\left(\frac{1\cdot 1}{2\cdot 4}\right)^{2}\left(\frac{a-b}{a+b}\right)^{4}+\left(\frac{1\cdot 1\cdot 3}{2\cdot 4\cdot 6}\right)^{2}\left(\frac{a-b}{a+b}\right)^{6}+\cdots\right\}   \\\end{array} $} 
\end{equation*}\hfill$\Box$
\end{thm}

The first formula for $p$ is due to \textsc{Maclaurin} \cite{mac} in 1742 while the second is due to \textsc{Ivory} \cite{I} in 1796.  \textsc{Berndt} (\cite{B}, p. 147) points out that \textsc{Ivory}'s formula can be obtained from a form of \textsc{Landen}'s transformation to be found in \textsc{Erd\'elyi}'s compendium (\cite{Erd}, p.111, formula (5)), but we have not been able to encounter a direct elementary proof in the literature which transforms the first of the two formulas for $p$ into the second one, and which does not appeal to general transformation formulas. \textsc{Ivory}'s own proof cleverly finesses this difficulty by showing that both series are equal to an integral.
\section{Three Lemmas}
\begin{lemma} The following identity is valid for all integers $k$ and real numbers $n$.\begin{equation*}
{-n-\frac{1}{2}\choose k}=\dfrac{\dis\binom{2n+2k}{n+k}\binom{n+k}{k}}{\dis\binom{2n}{n}}\frac{(-1)^{k}}{4^{k}}
\end{equation*}
\begin{proof}
The following identities between binomial coefficients are well known:\begin{align*}
\dis\binom{n}{k}\binom{n+\frac{1}{2}}{k}    & =\binom{2n+1}{k}\binom{2n+1-k}{k}\frac{1}{4^{k}}  \\
  &\\
   \binom{-r}{k} &=(-1)^{k}\binom{r+k-1}{k}  
\end{align*}
For example, see \textsc{Knuth} (\cite{K} pp. 71 and 46, resp.)
Now solve the first identity for $\binom{n+\frac{1}{2}}{k}$, write $(-n-1)$ in place of $n$ everywhere, use the second identity to get rid of the negative sign inside the binomial coefficients, and finally multiply and divide by $n!(n+k)!$
\end{proof}
\end{lemma}

This next lemma is the novelty in our proof.

\begin{lemma}
The following relation is valid for all integers $n\geqslant 1:$
\begin{equation}
S_{n}:=\sum_{m=0}^{n}\frac{(-1)^{m}}{2^{m}}\binom{2m}{m}\binom{n}{m}=\begin{cases}
    \dis \frac{1}{4^{\frac{n}{2}}}\binom{n}{\frac{n}{2}} & \text{if $n$ is even }, \\
     0 & \text{otherwise}.
\end{cases}
\end{equation}
\end{lemma}
\begin{proof}
We construct, in \emph{two different ways}, a certain polynomial, $G_{n}(x)$,  whose coefficient of $x^{n}$ is the sum $S_{n}$ defined in (3)\footnote{The idea behind the following proof is found in \textsc{Polya} and \textsc{Szeg\"o}, \cite{PS}, Part 1, problem \#39.}.

Observe that the \emph{coefficient of} $x^{m}$ in $$\left(1-\frac{x}{2}\right)^{2m}=\sum_{k=0}^{2m}\frac{(-1)^{k}}{2^{k}}\binom{2m}{k}x^{k}$$is equal to (put $k=m$)$$\frac{(-1)^{m}}{2^{m}}\binom{2m}{m}. $$Therefore, \emph{the coefficient of} $$x^{n}=x^{m}\cdot x^{n-m}$$in$$\left(1-\frac{x}{2}\right)^{2m}\binom{n}{m}x^{n-m}$$is equal to$$\frac{(-1)^{m}}{2^{m}}\binom{2m}{m}\binom{n}{m}.$$Thus, \emph{the sum} we are seeking to evaluae, viz.$$S_{n}:=\sum_{m=0}^{n}\frac{(-1)^{m}}{2^{m}}\binom{2m}{m}\binom{n}{m}$$is the \emph{coefficient} of $x^{n}$ in the \emph{sum}$$G_{n}(x):=\sum_{m=0}^{n}\left(1-\frac{x}{2}\right)^{2m}\binom{n}{m}x^{n-m}.$$This is the \textbf{\emph{first form}} of the polynomial $G_{n}(x).$

But $G_{n}(x)$ can also be written$$G_{n}(x)=\left\{\left(1-\frac{x}{2}\right)^{2}+x\right\}^{n}=\left(1+\frac{x^{2}}{4}\right)^{n},$$which is the \textbf{\emph{second form}} of $G_{n}(x)$. Expanding it by the binomial theorem we obtain$$G_{n}(x)=\left(1+\frac{x^{2}}{4}\right)^{n}=\sum_{m=0}^{n}\frac{1}{4^{m}}\binom{n}{m}x^{2m}.$$We see that all the powers of $x$ are \emph{even}.  Thus, if $n$ is \emph{odd}, the coefficient of $x^{n}$ is \emph{zero}, while if it is \emph{even}, the coefficient of $$x^{n}=(x^{\frac{n}{2}})^{2}$$is equal to$$ \frac{1}{4^{\frac{n}{2}}}\binom{n}{\frac{n}{2}}.$$
\end{proof}

The following result is well known (\cite{K}):

\begin{lemma}(\textbf{Vandermonde's Theorem}) The following relation holds for all real numbers $a$ and $b$:$$
\sum_{m=0}^{n}\binom{a}{m}\binom{b}{n-m}=\binom{a+b}{n}.$$
\end{lemma}\hfill$\Box$

\section{Sketch of the Proof}

Our proof of \textbf{Theorem 1} is based on an idea of \textsc{Carlson} \cite{C}.  However we use the fundamental \textbf{Lemma 2}, of which we have given a direct elementary proof based on the elementary properties of the binomial coefficients (see \S2), while \textsc{Carlson} proves the transformation on the basis of general expansion theorems for ``R-functions."  

Our proof consists of two steps:
\subsubsection*{Step 1}The left-hand side of \textbf{Theorem 1} is rearranged into the following power series in $\dis\left(\frac{2\sqrt{x}}{1+\sqrt{x}}\right)$:\begin{equation}
(1+x)^{\frac{1}{2}}\sum_{n=0}^{\infty}\frac{-1}{4n-1}\frac{1}{4^{2n}}\binom{4n}{2n}\frac{1}{4^{n}}\binom{2n}{n}\left(\frac{2\sqrt{x}}{1+\sqrt{x}}\right)^{2n}.
\end{equation}We achieve this via the binomial theorem, \textbf{Lemma 1} and \textbf{Lemma 2}.
\subsubsection*{Step 2}The power series $(2)$ is rearranged into a power series in $x$ alone via the binomial theorem.  Applying \textbf{Lemma 1} and \textsc{Vandermonde}'s theorem, the coefficients of this series collapse down to the coefficients of the right hand side of \textbf{Theorem 1}.

\section{Step 1}

All subsequent series converge absolutely and uniformly for $0\leqslant x\leqslant 1$, so that the following series manipulations are valid in that range.\begin{align*}
    &\dis(1+\sqrt{x})\left\{1-\frac{1}{2^{2}}\left(\frac{2\sqrt[4]{x}}{1+\sqrt{x}}\right)^{2}-\frac{1^{2}\cdot 3}{2^{2}\cdot 4^{2}}\left(\frac{2\sqrt[4]{x}}{1+\sqrt{x}}\right)^{4}-\frac{1^{2}\cdot 3^{2}\cdot 5}{2^{2}\cdot 4^{2}\cdot 6^{2}}\left(\frac{2\sqrt[4]{x}}{1+\sqrt{x}}\right)^{6}-\cdots\right\}\\
    &=\dis(1+\sqrt{x})\sum_{m=0}^{\infty}\frac{-1}{2m-1}\left\{\frac{1}{4^{m}}\binom{2m}{m}\right\}^{2}\left\{\dis\left(\frac{2\sqrt[4]{x}}{1+\sqrt{x}}\right)^{2}\right\}^{m}\\
&=\dis(1+\sqrt{x})\sum_{m=0}^{\infty}\frac{-1}{2m-1}\left\{\frac{1}{4^{m}}\binom{2m}{m}\right\}^{2}\dis\left(\frac{2\sqrt[4]{x}}{1+\sqrt{x}}\right)^{2m}   \\
    &=\dis\sum_{m=0}^{\infty}\frac{-1}{2m-1}\left\{\frac{1}{4^{m}}\binom{2m}{m}\right\}^{2}4^{m}(\sqrt{x})^{m}(1+\sqrt{x})^{-2m+1} \\ 
&=\dis\sum_{m=0}^{\infty}\frac{-1}{2m-1}\left\{\frac{1}{4^{m}}\binom{2m}{m}\right\}^{2}4^{m}(\sqrt{x})^{m}(1+2\sqrt{x}+x)^{-m+\frac{1}{2}}\\
&=\dis\sum_{m=0}^{\infty}\frac{-1}{2m-1}\left\{\frac{1}{4^{m}}\binom{2m}{m}\right\}^{2}4^{m}(\sqrt{x})^{m}(1+x)^{\frac{1}{2}}(1+x)^{-m}\left(1+\frac{2\sqrt{x}}{1+x}\right)^{-m+\frac{1}{2}}\\
&=-\sum_{m=0}^{\infty}\frac{1}{2m-1}\left\{\frac{1}{4^{m}}\binom{2m}{m}\right\}^{2}4^{m}\left(\frac{\sqrt{x}}{1+x}\right)^{m}(1+x)^{\frac{1}{2}}\sum_{k=0}^{\infty}\binom{-m+\frac{1}{2}}{k}\left(\frac{2\sqrt{x}}{1+x}\right)^{k}\\
&=-(1+x)^{\frac{1}{2}}\sum_{m=0}^{\infty}\frac{1}{2m-1}\left\{\frac{1}{4^{m}}\binom{2m}{m}\right\}^{2}2^{m}\sum_{k=0}^{\infty}\binom{-m+\frac{1}{2}}{k}\left(\frac{2\sqrt{x}}{1+x}\right)^{k+m}\\
&=-(1+x)^{\frac{1}{2}}\sum_{m=0}^{\infty}\sum_{k=0}^{\infty}\frac{2^{m}}{2m-1}\left\{\frac{1}{4^{m}}\binom{2m}{m}\right\}^{2}\binom{-m+\frac{1}{2}}{k}\left(\frac{2\sqrt{x}}{1+x}\right)^{k+m}\\
&=-(1+x)^{\frac{1}{2}}\sum_{m=0}^{\infty}\sum_{k=0}^{\infty}\frac{2^{m}}{2m-1}\left\{\frac{1}{4^{m}}\binom{2m}{m}\right\}^{2}\frac{2m-1}{2m-1+2k}\binom{-m-\frac{1}{2}}{k}\left(\frac{2\sqrt{x}}{1+x}\right)^{k+m}\\
&\left(\text{since $\binom{r}{k}=\dfrac{r}{r-k}\binom{r-1}{k}$}\right)\\
&=-(1+x)^{\frac{1}{2}}\sum_{m=0}^{\infty}\sum_{k=0}^{\infty}\frac{2^{m}}{4^{2m}}\binom{2m}{m}^{2}\frac{1}{2m+2k-1}\frac{\dis\binom{2m+2k}{m+k}\binom{m+k}{k}}{\dis\binom{2m}{m}}\frac{(-1)^{k}}{4^{k}}\left(\frac{2\sqrt{x}}{1+x}\right)^{k+m}\\
&(\text{by \textbf{Lemma 1}})\\
&=(1+x)^{\frac{1}{2}}\sum_{m=0}^{\infty}\sum_{k=0}^{\infty}\frac{2^{m}}{4^{2m}}\binom{2m}{m}\frac{1}{2m+2k-1}\dis\binom{2m+2k}{m+k}\binom{m+k}{k}\dis\binom{2m}{m}\frac{(-1)^{k}}{4^{k}}\left(\frac{2\sqrt{x}}{1+x}\right)^{k+m}\\
&(\text{now put $k:=n-m$ and interchange the order of summation})
\end{align*}\newpage\begin{multline*}
=(1+x)^{\frac{1}{2}}\sum_{n=0}^{\infty}\sum_{m=0}^{\infty}\frac{2^{m}}{4^{2m}}\binom{2m}{m}\frac{1}{2m+2(n-m)-1}\dis\binom{2m+2n-2m}{m+n-m}\binom{m+n-m}{n-m}\dis\binom{2m}{m}\times\\
    \times \frac{(-1)^{n-m}}{4^{n-m}}\left(\frac{2\sqrt{x}}{1+x}\right)^{n} 
\end{multline*}\begin{align*}
     &=(1+x)^{\frac{1}{2}}\sum_{n=0}^{\infty}\frac{(-1)^{n-1}}{2n-1}\frac{1}{4^{n}}\binom{2n}{n}\left\{\sum_{m=0}^{\infty}\frac{(-1)^{m}}{2^{m}}\binom{2m}{m}\binom{n}{n-m}\right\}\left(\frac{2\sqrt{x}}{1+x}\right)^{n}   \\
    &=(1+x)^{\frac{1}{2}}\sum_{n=0}^{\infty}\frac{(-1)^{n-1}}{2n-1}\frac{1}{4^{n}}\binom{2n}{n}\left[\begin{cases}
    \dis \frac{1}{4^{\frac{n}{2}}}\binom{n}{\frac{n}{2}} & \text{if $n$ is even }, \\
     0 & \text{otherwise}
\end{cases}
\right] \left(\frac{2\sqrt{x}}{1+x}\right)^{n} (\text{by \textbf{Lemma 2}})\\
 &=(1+x)^{\frac{1}{2}}\sum_{m=0}^{\infty}\frac{-1}{4m-1}\frac{1}{4^{2m}}\binom{4m}{2m}\dis \frac{1}{4^{m}}\binom{2m}{m}\left(\frac{2\sqrt{x}}{1+x}\right)^{2m} (\text{putting $n=2m$}.)
\end{align*}

Putting together the two extreme members of this chain of equalities, \textbf{\emph{we have shown:}}\begin{thm}\begin{align*}
 &\dis(1+\sqrt{x})\left\{1-\frac{1}{2^{2}}\left(\frac{2\sqrt[4]{x}}{1+\sqrt{x}}\right)^{2}-\frac{1^{2}\cdot 3}{2^{2}\cdot 4^{2}}\left(\frac{2\sqrt[4]{x}}{1+\sqrt{x}}\right)^{4}-\frac{1^{2}\cdot 3^{2}\cdot 5}{2^{2}\cdot 4^{2}\cdot 6^{2}}\left(\frac{2\sqrt[4]{x}}{1+\sqrt{x}}\right)^{6}-\cdots\right\}\\   
    &=(1+x)^{\frac{1}{2}}\sum_{m=0}^{\infty}\frac{-1}{4m-1}\frac{1}{4^{2m}}\binom{4m}{2m}\dis \frac{1}{4^{2m}}\binom{2m}{m}\left(\frac{2\sqrt{x}}{1+x}\right)^{2m}.  
\end{align*}
\end{thm}\hfill$\Box$

This completes the proof of \textbf{Step 1}.
\section{Step 2}

We expand the term $(1+x)^{-2m+\frac{1}{2}}$ by the binomial theorem.  Therefore, the RHS of the equation in \textbf{Theorem 3} becomes\begin{multline}
\sum_{m=0}^{\infty}\frac{-1}{4m-1}\frac{1}{4^{2m}}\binom{4m}{2m}\dis \frac{1}{4^{m}}\binom{2m}{m}4^{m}x^{m}\sum_{k=0}^{\infty}\binom{-2m+\frac{1}{2}}{k}x^{k}   \\
    =\sum_{m=0}^{\infty}\sum_{k=0}^{\infty}\frac{-1}{4m-1}\frac{1}{4^{2m}}\binom{4m}{2m}\dis \binom{2m}{m}\binom{-2m+\frac{1}{2}}{k}x^{m+k}  
\end{multline}Again, we put $k:=n-m$, interchange the order of summation, and observe that \textbf{\emph{the coefficient of}} $x^{m+k}=x^{n}$ \textbf{\emph{is}}\begin{align*}
    &=\frac{-1}{4m-1}\frac{1}{4^{2m}}\binom{4m}{2m}\dis \binom{2m}{m}\binom{-2m+\frac{1}{2}}{k}   \\
    &=\frac{-1}{4m-1}\frac{1}{4^{2m}}\binom{4m}{2m}\dis \binom{2m}{m}\frac{4m-1}{4m+2k-1}\binom{-2m-\frac{1}{2}}{k}\\
  &=\frac{-1}{4^{2m}}\binom{4m}{2m}\dis \binom{2m}{m}\frac{(-1)^{k}}{4^{k}}\frac{\dis\binom{4m+2k}{m+k}\binom{2m+k}{k}}{\dis\binom{4m}{2m}}\frac{1}{4m+2k-1}\\
 &=\frac{-1}{4^{2m}}\binom{2m}{m}\frac{(-1)^{n-m}}{4^{n-m}}\binom{2m+2n}{m+n}\binom{n+m}{n-m}\frac{1}{2n+2m-1}\\ 
 &=\frac{(-1)^{n-m+1}}{4^{n+m}}\binom{2n+2m}{n+m}\frac{(n+m)!}{(n-m)!(2m)!}\frac{(2m)!}{(m)!(m)!}\frac{1}{2n+2m-1}\\
&=\frac{(-1)^{n-m+1}}{4^{n+m}}\binom{2n+2m}{n+m}(n+m)!\frac{1}{m!(n-m)!}\frac{1}{m!}\frac{1}{2n+2m-1}\\
 &=\frac{(-1)^{n-m+1}}{4^{n+m}}\binom{2n+2m}{n+m}(n+m)!\frac{1}{m!(n-m)!}\frac{(n-m)!}{m!}\binom{n}{n-m}\frac{1}{2n+2m-1}\\ 
 &=\frac{1}{2n-1}\frac{(-1)^{n-m+1}}{4^{n}}\binom{2n}{n}\dfrac{\dis\frac{1}{4^{m}}\binom{2n+2m}{n+m}\frac{(n+m)!}{n!m!}\binom{n}{n-m}}{\dis\frac{1}{2n-1}\binom{2n}{n}}\frac{1}{2n+2m-1}\\ 
&=\frac{(-1)^{n-1}}{2n-1}\frac{1}{4^{n}}\binom{2n}{n}\binom{-n+\frac{1}{2}}{m}\binom{n}{n-m},
\end{align*}where we have used both the conclusion and the proof of \textbf{Lemma 1}.  So, if we now substitute this last expression:$$\frac{(-1)^{n-1}}{2n-1}\frac{1}{4^{n}}\binom{2n}{n}\binom{-n+\frac{1}{2}}{m}\binom{n}{n-m}$$into the RHS of (5) we obtain\begin{align*}
&\dis(1+\sqrt{x})\left\{1-\frac{1}{2^{2}}\left(\frac{2\sqrt[4]{x}}{1+\sqrt{x}}\right)^{2}-\frac{1^{2}\cdot 3}{2^{2}\cdot 4^{2}}\left(\frac{2\sqrt[4]{x}}{1+\sqrt{x}}\right)^{4}-\frac{1^{2}\cdot 3^{2}\cdot 5}{2^{2}\cdot 4^{2}\cdot 6^{2}}\left(\frac{2\sqrt[4]{x}}{1+\sqrt{x}}\right)^{6}-\cdots\right\}\\
    &=\sum_{n=0}^{\infty}\sum_{m=0}^{n}\frac{(-1)^{n-1}}{2n-1}\frac{1}{4^{n}}\binom{2n}{n}\binom{-n+\frac{1}{2}}{m}\binom{n}{n-m}x^{n}\\
&=\sum_{n=0}^{\infty}\frac{(-1)^{n-1}}{2n-1}\frac{1}{4^{n}}\binom{2n}{n}\left\{\sum_{m=0}^{n}\binom{-n+\frac{1}{2}}{m}\binom{n}{n-m}\right\}x^{n}\\
&(\text{and by \textsc{Vandermonde}'s theorem this is equal to})      
\end{align*}\newpage\begin{align*}
  &=\sum_{n=0}^{\infty}\frac{(-1)^{n-1}}{2n-1}\frac{1}{4^{n}}\binom{2n}{n}\binom{\frac{1}{2}}{n}x^{n}   \\
    &=\sum_{n=0}^{\infty}\frac{(-1)^{n-1}}{2n-1}\frac{1}{4^{n}}\binom{2n}{n}\frac{1}{4^{n}}\frac{(-1)^{n-1}}{2n-1}\binom{2n}{n}x^{n}\\ 
    &=\sum_{n=0}^{\infty}\left\{\frac{(-1)^{n-1}}{2n-1}\frac{1}{4^{n}}\binom{2n}{n}\right\}^{2}x^{n} \\
    &=\dis1+\left(\frac{1}{2}\right)^{2}x+\left(\frac{1\cdot 1}{2\cdot 4}\right)^{2}x^{2}+\left(\frac{1\cdot 1\cdot 3}{2\cdot 4\cdot 6}\right)^{2}x^{3}+\left(\frac{1\cdot 1\cdot 3\cdot 5}{2\cdot 4\cdot 6\cdot 8}\right)^{2}x^{4}+\cdots
\end{align*}which is the RHS of \textbf{Theorem 1}.  This completes the proof of \textbf{Theorem 1}.\hfill$\Box$


\end{document}